\documentclass {amsart}

\usepackage{mathrsfs}
\usepackage{graphicx,verbatim}
\usepackage{amscd}

\usepackage{amsmath}
\usepackage{amssymb}

\newcommand{\bbP}{\mathbb{P}}
\newcommand{\buk}{\noindent}
\newcommand{\bbF}{\mathbb{F}}
\newcommand{\bbZ}{\mathbb{Z}}
\newcommand{\hf}{\hfill$\Box$}
\newcommand{\calM}{\mathcal{M}}
\newcommand{\calS}{\mathcal{S}}

\newcommand{\calO}{\mathcal{O}}
\newcommand{\Prf}{\buk{\bf Proof}}

\newtheorem{theorem}{Theorem}
\newtheorem{lemma}[theorem]{Lemma}
\newtheorem{proposition}[theorem]{Proposition}
\newtheorem{definition}[theorem]{Definition}
\newtheorem{corollary}[theorem]{Corollary}
\newtheorem{conclusion}[theorem]{Conclusion}
\newtheorem{remark}[theorem]{Remark}
\newtheorem{example}[theorem]{Example}
\newtheorem{notation}[theorem]{Notation}
\newtheorem{convention}[theorem]{Convention}

\begin{document}


\title[Moduli of bundles I]{Moduli of bundles over rational surfaces and elliptic curves
I: simply laced cases }


\author{Naichung Conan Leung}


\address{the Institute of Mathematical Sciences and Department of Mathematics,
         The Chinese University of Hong Kong,
         Shatin, N.T., Hong Kong}

\email{leung@ims.cuhk.edu.hk}
\thanks{The first author is partially supported by Hong Kong RGC Grant. The second author is supported
by the SFB/TR 45: Periods, Moduli Spaces and Arithmetic of Algebraic
Varieties of the DFG.}


\author{Jiajin Zhang}
\address{Department of Mathematics, Sichuan University, Chengdu, 610000, P.R. China}
\email{jjzhang@scu.edu.cn}
\curraddr{Institute of Mathematics, Johannes Gutenberg-University Mainz, Mainz, 55099, Germany}
\email{zhangji@uni-mainz.de}


\date{Nov. 20, 2008.}



\begin{abstract}
  It is well-known that del Pezzo surfaces of
degree $9-n$ one-to-one correspond to flat $E_n$ bundles over an
elliptic curve. In this paper, we construct $ADE$ bundles over a
broader class of rational surfaces which we call $ADE$ surfaces, and
extend the above correspondence to all flat $G$ bundles over an
elliptic curve, where $G$ is any simply laced, simple, compact and
simply-connected Lie group. In the sequel, we will construct $G$
bundles for non-simply laced Lie group $G$ over these rational
surfaces, and extend the above correspondence to non-simply laced
cases.
\end{abstract}

\maketitle

\section*{Introduction}

        Let $S$ be a smooth rational surfaces. If
     the anti-canonical line bundle  $-K_S$ is ample, then $S$ is called a
     {\it del Pezzo surface}. It is well-known that a del Pezzo surface can be
     classified as a blow-up of $\mathbb{CP}^2$ at $n(n\leq 8)$
     points in general position or $\mathbb{CP}^1\times \mathbb{CP}^1$. When these blown-up points are in
     {\it almost general position}, such a surface is called a {\it generalized del Pezzo surface},
     according to Demazure \cite{Demazure}. It is
     also well-known that the sub-lattice $K_S^{\perp}$ of $Pic(S)$
     is a root lattice of type $E_n$. For more results on (generalized) del Pezzo surfaces one can
     see \cite{Demazure} and \cite{Manin}. Thus there is a natural Lie algebra bundle of type $E_n$ over $S$.
     By restriction to a fixed smooth anti-canonical curve $\Sigma$, one obtains a flat $E_n$ bundle over $\Sigma$.
     Moreover, Donagi \cite{Donagi} \cite{Donagi2} and Friedman-Morgan-Witten \cite{FMW1}  \cite{FMW2} prove that
     the moduli space of del Pezzo surfaces with fixed anti-canonical curve $\Sigma$ can be identified with the moduli
     space of  flat $E_n$ bundles over this elliptic curve $\Sigma$.

      In this paper, we will extend this correspondence to all
     compact, simple, simply laced and simply connected Lie groups and to a broader class
     of rational surfaces, which are called $ADE$ surfaces. This paper contains parts of the preprint \cite{Leung},
     especially the construction of Lie algebra bundles and their (fundamental) representation bundles,
     and we shall refer to \cite{Leung} for some of
     the proofs. Next we sketch the contents briefly. \\

        In Section 1, we first analyze the structure  of the Picard lattice of a
        rational surface which is a blow-up of $\bbP^2$, $\bbP^1\times \bbP^1$ or
    the Hirzebruch surface $\bbF_1$ at some points. We shall see that
    there is a sub-lattice of the Picard lattice which is a root lattice of
    $ADE$-type.

        Next we generalize the definition of del Pezzo surfaces to that of $ADE$ surfaces,
        where an $E_n$ surface is just a del Pezzo surface of degree $9-n$.
        Roughly speaking, an $ADE$ surface  $S$ is a rational surface with a smooth rational curve $C$ on
     $S$ such that the sub-lattice $\langle K_S,C\rangle^{\perp}$ of $Pic(S)$ is
     an irreducible root lattice (see Definition~\ref{ADE-surface}).  The condition in Definition~\ref{ADE-surface}
     implies that $C^2=-1,0$ or $1$, and that the sub-lattice $\langle K_S,C\rangle^{\perp}$ is a root lattice of type
     $E_n$, $D_n$, or $A_n$ respectively (Proposition~\ref{ADE-class}).
        Therefore such a surface is called a rational surface of  $E_n$-type, $D_n$-type, or $A_n$-type accordingly.

          Note that the definition of  an $E_n$ surface implies that after blowing down the $(-1)$ curve $C$,
        the anti-canonical line bundle $-K$ will be ample. So the
        resulting surface is just a del Pezzo surface.
        Thus the definition of $ADE$ surfaces naturally generalizes that of
        del Pezzo surfaces.\par

       After this, we prove that an $ADE$ surface is nothing but a blow-up of $\bbP^2$, $\bbP^1\times \bbP^1$ or
    $\bbF_1$ at some points in general position. This gives us an explicit construction for any $ADE$
    surface.\\

        In Section 2, we construct Lie algebra bundles of $ADE$-type, and their natural representation bundles over
     those surfaces discussed in Section 1.
        By a Lie algebra bundle over a surface $S$, we mean a vector bundle
     which has a fiberwise  Lie algebra structure, and this
     structure is compatible with any trivialization. Similarly, by a
     representation bundle, we mean a vector bundle which is a
     fiberwise representation of a Lie algebra bundle, and this fiberwise representation is compatible with any
     trivialization.\par
       More precisely, let  $S$ be an $ADE$ surface. Since the sub-lattice $\langle K_S,C\rangle^{\perp}$ of $Pic(S)$ is a
     root lattice, we can explicitly construct a natural Lie algebra bundle of
     corresponding type over $S$, using the root system of
     the root lattice $\langle K_S,C\rangle^{\perp}$.  Using the lines and rulings on $S$, we can also construct
     natural fundamental representation bundles over $S$.\\

        In Section 3, we relate the above Lie algebra bundles of $ADE$-type over $ADE$ rational surfaces
    to flat $G$ bundles over an elliptic curve $\Sigma$,
    where $G$ is a compact Lie group of corresponding type. If an $ADE$ rational surface $S$ contains
    a fixed smooth elliptic curve $\Sigma$ as an
   anti-canonical curve, then by restriction, one obtains flat $ADE$-bundles
   over $\Sigma$. We can prove this restriction identifies the
   moduli space of flat $ADE$ bundles over $\Sigma$ and the moduli
   space of the pairs $(S,\Sigma\in |-K_S|)$ with extra structure $\zeta_G$ which
   is called a $G$-{\it configuration} (Definition~\ref{config}).  Our main result in this paper is the
   following theorem.

   \begin{theorem} Let $\Sigma$ be a fixed elliptic curve,  and let $G$ be a simple, compact, simply laced and simply
   connected Lie group. Denote $\calS(\Sigma,G)$
   the moduli space of the pairs $(S,\Sigma)$, where $S$ is an $ADE$ rational surface with $\Sigma\in |-K_S|$. Denote
   $\calM_\Sigma^G$ the moduli space of flat $G$-bundles over $\Sigma$. Then by restriction, we have\par
     (i) $\calS(\Sigma,G)$ can be  embedded into $\calM_\Sigma^G$  as an open dense subset.\par
     (ii) There exists a natural and explicit compactification for $\calS(\Sigma,G)$, denoted by
     $\overline{\calS(\Sigma,G)}$,  such that this embedding
     can be extended to an isomorphism from $\overline{\calS(\Sigma,G)}$ onto $\calM_\Sigma^G$.\par
   (iii) Any surface corresponding to a boundary point in $\overline{\calS(\Sigma,G)}\, \backslash \calS(\Sigma,G)$
   is equip-ped with a $G$-configuration, and on such a surface,
   any smooth rational curve has a self-intersection number at
   least $-2$. Furthermore, in $E_n$ case, all $(-2)$ curves form  chains of
   $ADE$-type, and the anti-canonical model of such a surface admits at worst $ADE$-singularities.
   \end{theorem}

      Physically, when $G=E_n$ is a simple subgroup of $E_8\times E_8$, these $G$ bundles are related
   to the duality between  $F$-theory and string theory. Among other things, this duality predicts the moduli of
   flat $E_n$ bundles over a fixed elliptic curve $\Sigma$ can be identified with the
   moduli of del Pezzo surfaces with fixed anti-canonical curve
   $\Sigma$. For details, one can consult  \cite{Donagi}
   \cite{Donagi2}  \cite{FMW1} and  \cite{FMW2}. Our result can be considered as a test of above duality for
    other Lie groups.

     As an application, we have a more intuitive explanation for the well-known moduli space $\calM_\Sigma^G$ of
     flat $G$-bundles over a fixed elliptic curve $\Sigma$. And we can see very clearly how the Weyl group
     of $G$ acts on the marked moduli space of flat $G$-bundles over $\Sigma$.

\begin{notation} {\rm In this paper, we will fix some notations from Lie theory. Let
$G$ be a compact, simple and simply-connected Lie group. We
denote\par
   $r(G)$: the rank of $G$;\par
   $R(G)$: the root system;\par
   $R_c(G)$: the coroot system;\par
   $W(G)$: the Weyl group;\par
   $\Lambda(G)$: the root lattice;\par
   $\Lambda_c(G)$: the coroot lattice;\par
   $\Lambda_w(G)$: the weight lattice;\par
   $T(G)$: a maximal torus;\par
   $ad(G)$: the adjoint group of $G$, i.e. $G/C(G)$ where $C(G)$
   is the center of $G$;\par
   $\Delta(G)$: a simple root system of $G$.\par
   When there is no confusion, we just ignore  the letter $G$.}
   \end{notation}

\noindent{\bf Acknowledgements}. The revision of this paper was done
during the second author's stay at the Institute of mathematics of
Johannes Gutenberg-University Mainz. The second author would like to
thank the Institute (especially Prof. Kang Zuo) for their kind
hospitality and good working conditions. He would also like to
express his gratitude to Dr. Changzheng Li and Prof. Xiaowei Wang
for many discussions. We would like to thank the referee for showing
us the relation (see Remark~\ref{deformation}) of our work with
Looijenga and Pinkham's work on the semi-universal deformation of a
simple-elliptic singularity
\cite{Looijenga1}\cite{Looijenga2}\cite{Looijenga98}\cite{Pinkham},
and for useful comments that improved the exposition of this paper
significantly.

\section{Rational surfaces of $ADE$-type}

       Before defining what $ADE$ surfaces are, we first give their
   explicit constructions.

\subsection{}
      First consider the $E_n$ case, that is, the case of del Pezzo surfaces.
   We start with a complex projective plane $\mathbb{P}^2$ and $n$ points $x_1,\cdots,x_n$ on $\bbP^2$ with $n\leq
   8$. Note that $x_2,\cdots,x_n$ may be {\it infinitely near} points. For example, we say that $x_2$ is  {\it infinitely near} $x_1$
   if $x_2$ lies on the exceptional curve obtained by blowing up $x_1$.
   Blowing up $\mathbb{P}^2$ at these points in turn, we obtain a rational surface, denoted
   $X_n(x_1,\cdots,x_n)$ or $X_n$ for brevity.

        These points are said to be {\it in general position} if they satisfy the following
    conditions:\par
    (i) They are distinct points;\par
    (ii) No three of them are collinear;\par
    (iii) No six of them lie on a common conic curve;\par
    (iv) No cubics pass through $8$ points with one of them a double point.  \par
    The following result is well-known (see \cite{Demazure} and \cite{Manin}).

    \begin{lemma}\label{equiv-del-pezzo}  Let $x_i\in \mathbb{P}^2, i=1,\cdots,n,n\leq 8$.  Then the following
    conditions are equivalent:\par
     (i)  These points are in general position.\par
    (ii)  The self-intersection number of any rational curve on $X_n$ is bigger than or equal to $-1$.\par
    (iii) The anti-canonical class $-K_{X_n}$ is
    ample.$\hfill\Box$\end{lemma}

    A surface $X_n$ is called a {\it del Pezzo surface} if it
    satisfies one of the above equivalent conditions.\\

    We say that $x_i\in \mathbb{P}^2, i=1,\cdots,n$ with $n\leq 8$ are
    in {\it almost general position} if any smooth rational curve on $X_n$ has a self-intersection number at least
    $-2$,  and such a surface is called a {\it generalized del Pezzo surface} (see \cite{Demazure}).\\

     Let $h$ be the class of lines in $\mathbb{P}^2$ and  $l_i$ be
   the exceptional divisor corresponding to
   the blow-up at  $x_i\in \mathbb{P}^2,i=1,\cdots, n$. Denote $Pic(X_n)$ the Picard group of $X_n$, which is
   isomorphic to $H^2(X_n, \mathbb{Z})$.  Then  $Pic(X_n)$ is a lattice  with basis
   $h,l_1,\cdots,l_n$, of signature $(1,n)$.  Let $K=-3h+l_1+\cdots+l_n$ be the canonical class.
   We extend the definition of the  Lie algebras $E_n,n=6,7,8$ to all $n$ with $0\leq n\leq 8$ by
   setting $E_0=0, E_1=\mathbb{C}, E_2=A_1\times \mathbb{C}, E_3=A_1\times A_2,E_4=A_4$  and $E_5=D_5$.\par

   Denote
   \begin{eqnarray}P_n &=& \{x\in Pic(X_n)\ |\  x\cdot K=0\},\nonumber
   \\ R_n &=& \{x\in  Pic(X_n)\ |\ x\cdot K=0,\ x^2=-2\}\subset P_n,\nonumber\\
   I_n &=& \{x\in  Pic(X_n)\ |\  x^2=-1=x\cdot K\},\mbox{ and }\nonumber\\
   C_n &=& \{\zeta_n=(e_1,\cdots,e_n)\ |\  e_i\in I_n,\ e_i\cdot e_j=0,\ i\neq j\}.\nonumber
   \end{eqnarray}

      An element of $I_n$ is called an {\it exceptional divisor}, and an
      element $\zeta_n\in C_n$ is called an {\it exceptional system} (of divisors) (see \cite{Demazure} and
      \cite{Manin}).

   \begin{lemma} \label{root-E} (i) $R_n$ is a root system of type $E_n$ with a system of simple roots
   $\alpha_1=l_1-l_2,\ \alpha_2=l_2-l_3,\ \alpha_3=h-l_1-l_2-l_3,\ \alpha_4=l_3-l_4,\ \cdots,\ \alpha_n=l_{n-1}-l_n $.
   Its root lattice is just   $P_n$, and its weight lattice is $Q_n=Pic(X_n)/\mathbb{Z}K$.
   Let $l\in I_n$, then $R_n\cap l^{\perp}$ is
   a root system of type $E_{n-1}$, and $P_n\cap l^{\perp}$ is its root
   lattice.\par (ii) The Weyl group $W(E_n)$ acts on $C_n$ simply
   transitively.\end{lemma}

   \noindent{\bf Proof}. (i) For the proof that $R_n$ is a root system of type $E_n$ with given simple roots,
   see Manin's book \cite{Manin}.   $Pic(X_n)$ is a lattice with $\mathbb{Z}$-basis $h,l_1,\cdots,l_n$.
   Obviously, $\{e_0=l_1,e_1=\alpha_1,\cdots,e_n=\alpha_n\}$
    forms another $\mathbb{Z}$-basis. Take any $x\in P_n\subset Pic(X_n)$.
    Let $x=\sum a_i\cdot e_i$. Then $x\cdot K=0$ implies $a_0=0$. So $P_n$ is the  root lattice of $R_n$. \par
    The natural pairing $P_n\otimes Pic(X_n)\rightarrow\mathbb{Z}$ induces a perfect
    pairing
    $$P_n\otimes (Pic(X_n)/\mathbb{Z}K)\rightarrow\mathbb{Z}.$$
     So the weight lattice is just $Pic(X_n)/\mathbb{Z}K$.\par
       For the last assertion, we can assume $l=l_8$, then it is true obviously.  \par
    (ii) See \cite{Manin}. \hfill$\Box$\\

     The Dynkin diagram is the following\\
\\
   \begin{center}
   \includegraphics{ade.2}\\
   \end{center}

\subsection{}

      Next we consider the $D_n$ case. Let $Y=\mathbb{F}_1$ be a {\it Hirzebruch surface}, and
   fix the ruling $f$ and the section $s$, where $s^2=-1$. In fact
   $\mathbb{F}_1$ is the blow-up of $\bbP^2$ at one point $x_0$. Thus $f=h-l_0,\ s=l_0$ where $h$ is
   the class of lines on $\mathbb{P}^2$ and  $l_0$ is the exceptional curve. Blowing up $Y$ at $n$ points $x_1,\cdots,x_n$ we
   obtain $Y_n$. The Picard group $Pic(Y_n)$ of $Y_n$ is $H^2(Y_n, \mathbb{Z})$,
   which  is a lattice with basis
   $s,f,l_1,\cdots,l_n$. The canonical class $K=-(2s+3f-\sum\limits_{i=1}^{n}l_i)$.\par
   Denote
   \begin{eqnarray}P_n &=& \{x\in Pic(Y_n)\ |\  x\cdot K=0=x\cdot f\},\nonumber
   \\ R_n &=& \{x\in Pic(Y_n)\ |\ x\cdot K=0=x\cdot f,\ x^2=-2\},\nonumber \\
   I_n &=& \{x\in Pic(Y_n)\ |\  x^2=-1=x\cdot K,\ x\cdot f=0\},\nonumber\\
   C_n &=& \{\zeta_n=(e_1,\cdots,e_n)\ |\  e_i\in I_n, e_i\cdot e_j=0, i\neq
   j,\nonumber\\
   & & \sum e_i\cdot s\equiv 0 \ mod \ 2\}.\nonumber
   \end{eqnarray}

        Similarly as before,  an element $\zeta_n\in C_n$ is called an {\it exceptional system} (of
   divisors).

   \begin{lemma} \label{root-D}  (i) $R_n$ is a root system of type $D_n$ with a system of simple roots
   $\alpha_1=f-l_1-l_2,\alpha_2=l_1-l_2,\cdots,\alpha_n=l_{n-1}-l_n $. Its root lattice is just
   $P_n$ and its weight lattice is $Q_n=Pic(Y_n)/\mathbb{Z}\langle f,K\rangle$.\par (ii) The Weyl group $W(D_n)$ acts on $C_n$ simply transitively.\end{lemma}

   \noindent{\bf Proof}. (i) $Pic(Y_n)$ is a lattice with $\mathbb{Z}$-basis
   $s,f,l_i,i=1,\cdots,n$. Let $x=as+bf+\sum c_i l_i\in R_n$ where $a,b,c_i\in\mathbb{Z}$. Then
   we have a system of linear equations
    $$\left\{
    \begin{array}{l}  x^2=-2, \\
     x\cdot K=0=x\cdot f.
    \end{array}\right.
    $$
    Solving this, we obtain
     $$\left\{
    \begin{array}{l}  a=0,\\ \sum c_i^2=2, \\ 2b=-\sum c_i.
    \end{array}\right.
    $$
         So, $x=\pm(l_i-l_j),i\neq j$ or $x=\pm(f-l_i-l_j),i\neq j$.
     That is $R_n=\{\pm(l_i-l_j), \pm(f-l_i-l_j)|\ i\neq j\}$.
     This implies that $R_n$ is a root system of $D_n$-type with
     indicated simple roots.

         Obviously, $\{e_1=s,e_2=l_1,e_{i+2}=\alpha_i,i=1,\cdots,n\}$
    forms another $\mathbb{Z}$-basis. Take any $x\in P_n\subset Pic(Y_n)$.
    Let $x=\sum a_i\cdot e_i$. Then $x\cdot K=0=x\cdot f$ implies $a_1=a_2=0$. So $P_n$ is the  root lattice of
    $R_n$.\par
         The natural pairing $P_n\otimes Pic(Y_n)\rightarrow\mathbb{Z}$ has kernel
    $\mathbb{Z}\langle f,-2s+\sum l_i\rangle=\mathbb{Z}\langle f,K\rangle$. So the pairing induces a perfect
    pairing $P_n\otimes (Pic(Y_n)/\mathbb{Z}\langle f,K\rangle)\rightarrow\mathbb{Z}$.
    Hence the weight lattice is just $Pic(Y_n)/\mathbb{Z}\langle f,K\rangle$. \par
    (ii) A simple computation shows that $$I_n=\{l_i,f-l_i|i=1,\cdots,n\}.$$  Thus all the elements of $C_n$
    are of the form $\zeta_n=(u_1,\cdots,u_n)$ where the number
    of $u_i$'s, such that $u_i=f-l_k$ for some $k$, is even. Then by the structure of $W(D_n)$, the result is clear. \hfill$\Box$\\

   The Dynkin diagram is the following\\
   \\
   \begin{center}
   \includegraphics{ade.3}\\
   \end{center}

\subsection{}

   In the following we consider the $A_{n-1}$ case. For this, let $Z_n$ be just the same as $Y_n$.\par
   Denote
   \begin{eqnarray}P_{n-1} &=& \{x\in Pic(Z_n)\ |\  x\cdot K=x\cdot f=x\cdot s=0\},\nonumber
   \\ R_{n-1} &=& \{x\in Pic(Z_n)\ |\ x\cdot K=x\cdot f=x\cdot s=0,\ x^2=-2\},\nonumber\\
   I_{n-1} &=& \{x\in Pic(Z_n)\ |\  x^2=-1=x\cdot K,\ x\cdot f=0=x\cdot s\},\nonumber\\
   C_{n-1} &=& \{\zeta_n=(e_1,\cdots,e_n)\ |\  e_i\in I_{n-1}, e_i\cdot e_j=0, i\neq
   j\}.\nonumber
   \end{eqnarray}

    As before, an element of  $\zeta_n\in C_{n-1}$ is called an {\it exceptional system} (of
    divisors).

\begin{lemma} \label{root-A}   (i) $R_{n-1}$ is a root system of type
   $A_{n-1}$ with a system of simple roots   $\alpha_1=l_1-l_2,\cdots,\alpha_{n-1}=l_{n-1}-l_n $.
   Its root lattice is just
   $P_{n-1}$ and its weight lattice is $Pic(Z_n)/\mathbb{Z}\langle f,s,K\rangle$.
   \par (ii) The Weyl group $W(A_{n-1})$ acts on $C_{n-1}$ simply transitively. In fact, \\ $W(A_{n-1})$ acts as
   the permutation group of $l_1,\cdots,l_n$.\par
   (iii) Let $e$ be a $(-1)$ curve which does not meet $s$. Then there exist $i,j$ with $i\neq j$ such that
   $e=s+f-l_i-l_j$, and when $n\geq 4$, $\langle K,s,f,e\rangle^{\perp}$ is a reducible root
   lattice of type $A_1\times A_{n-3}$; when $n=3$, $\langle K,s,f,e\rangle^{\perp}$ is not a root
   lattice; when $n=2$, $\langle K,s,f,e\rangle^{\perp}$ is the same as $P_1$, which is of type
   $A_1$.\par
   (iv) Let $e_i,1\leq i\leq k,\ k\geq 2$ be $(-1)$ curves such that $s,e_i,1\leq i\leq k$ are
   disjoint pairwise. Then when $k\neq 3$, $\langle K,s,f,e_i,1\leq i\leq k\rangle^{\perp}$ is not a root
   lattice. When $k=3$, (a) if $e_1=s+f-l_{i_2}-l_{i_3},\ e_2=s+f-l_{i_1}-l_{i_3},\ e_3=s+f-l_{i_1}-l_{i_2}$
   then $\langle K,s,f,e_1,e_2,e_3\rangle^{\perp}$ is a root lattice of $A$-type; (b) otherwise, $\langle K,s,f,e_1,e_2,e_3\rangle^{\perp}$
   is not a root lattice. \end{lemma}

 \noindent  {\bf Proof}. (i) $Pic(Z_n)$ is a lattice with $\mathbb{Z}$-basis
   $s,f,l_i,i=1,\cdots,n$. A simple computation shows that
   $$R_{n-1}=\{l_i-l_j\ |\ i\neq j\}.$$ Then it is obviously a root system of
   type $A_{n-1}$ with given simple roots.\par
   Obviously, $\{e_1=s,e_2=f,e_3=l_1,e_{i+3}=\alpha_i,i=1,\cdots,n\}$
    forms another $\mathbb{Z}$-basis. Take any $x\in P_{n-1}\subset Pic(Z_n)$.
    Let $x=\sum a_i\cdot e_i$. Then $x\cdot K=x\cdot f=x\cdot s=0$ implies $a_1=a_2=a_3=0$.
    So $P_{n-1}$ is the  root lattice of $R_{n-1}$. \par
    The natural pairing $$P_{n-1}\otimes Pic(Z_n)\rightarrow\mathbb{Z}$$ has a kernel
    $$\mathbb{Z}\langle f,s,\sum l_i\rangle=\mathbb{Z}\langle f,s,K\rangle.$$ So the pairing induces a perfect
    pairing $$P_{n-1}\otimes Pic(Z_n)/\mathbb{Z}\langle f,s,K\rangle)\rightarrow\mathbb{Z}.$$
    Hence the weight lattice is just $Pic(Z_n)/\mathbb{Z}\langle f,s,K\rangle$.\par
    (ii)  In fact  $I_{n-1}=\{l_1,\cdots,l_n\}$. So an element of $C_{n-1}$ is just a permutation of $l_1,\cdots,l_n$. \par
    (iii)  Let $e=as+bf+\sum c_il_i$, then $e$ is a $(-1)$ curve and $e\cdot
    s=0$ imply that $e$ must be of the form $s+f-l_i-l_j,i\neq j$.
    Without loss of generality, we can assume that $e=s+f-l_1-l_2$.
    Then the result follows from a simple computation. \par
    (iv) First let $k=2$. From the proof of (iii), we
    know both $e_1$  and $e_2$ are the form $s+f-l_i-l_j,i\neq
    j$. Since $e_1\cdot e_2=0$, we can assume
    $e_1=s+f-l_1-l_2$ and  $e_2=s+f-l_1-l_3$. Then the result
    follows easily.  For $k=3$, if $e_1=s+f-l_{i_2}-l_{i_3},\ e_2=s+f-l_{i_1}-l_{i_3},\ e_3=s+f-l_{i_1}-l_{i_2}$
   then $\langle K,s,f,e_1,e_2,e_3\rangle^{\perp}=\langle K,s,f,l_{i_1},l_{i_2},l_{i_3}\rangle^{\perp}$. We can assume
   $l_{i_1}=l_1, l_{i_2}=l_2, l_{i_3}=l_3$. Then  $\langle K,s,f,l_1,l_2,l_3\rangle^{\perp}$ is a root lattice of $A$-type.
   Other cases are similar. \hfill$\Box$\\

   The Dynkin diagram is the following\\
\\
   \begin{center}
   \includegraphics{ade.4}\\
   \end{center}

  \vspace{0.5cm}    Note that  Lemma~\ref{root-D} and Lemma~\ref{root-A} (i) (ii) are still true if we replace
  $\mathbb{F}_1$  by any {\it Hirzebruch surface} $\bbF_k(k\geq 0)$.

\subsection{}
     Now we show that in a suitable sense, the converse of the above lemmas is also true.
     As promised in the introduction, we will see that the following definition
   generalizes that of {\it del Pezzo surfaces}.

   \begin{definition} \label{ADE-surface} Let $(S,C)$ be a pair consisting of a  smooth rational surface $S$
   and a smooth rational curve $C\subset S$ with $C^2\neq 4$.  The pair $(S,C)$ is called
   of $ADE$-type (or an $ADE$ surface) if it satisfies the
   following two conditions:\par
     (i) Any (smooth) rational curve on $S$ has a
   self-intersection number at least  $-1$;\par
     (ii) The  sub-lattice  $\langle K_S,C\rangle^{\perp}$ of $Pic(S)$ is an irreducible root lattice
   of rank equal to $rank(Pic(S))-2$.\end{definition}

   The following proposition shows that such surfaces can be classified into three types.

\begin{proposition} \label{ADE-class} Let $(S,C)$ be a rational surface of
$ADE$-type. Let $n=rank(Pic(S))-2$. Then $C^2\in\{-1,0,1\}$ and
\par
     (i) when $C^2=-1$,  $\langle K_S,C\rangle^{\perp}$ is of $E_n$-type, where $4\leq n\leq
     8$;\par
     (ii) when $C^2=0$,  $\langle K_S,C\rangle^{\perp}$ is of $D_n$-type, where $n\geq 3$;\par
     (iii) when $C^2=1$,  $\langle K_S,C\rangle^{\perp}$ is of $A_n$-type.\end{proposition}

\Prf. By the first condition in Definition~\ref{ADE-surface},
$C^2\geq -1$. Therefore there are the following four cases.\par
Firstly, suppose $C^2=-1$. Then we can contract $C$ to obtain a
smooth surface $\widetilde{S}$. Let $\pi:S\rightarrow \widetilde{S}$
be the blow-down. Then the projection
$$Pic(S)=Pic(\widetilde{S})\oplus \mathbb{Z}\langle C\rangle\rightarrow
Pic(\widetilde{S})$$ induces an isomorphism $\langle
K_S,C\rangle^{\perp}\cong\langle K_{\widetilde{S}}\rangle^{\perp}$.
But the latter is an irreducible root system if and only if
$\widetilde{S}$ is a blow-up of $\mathbb{CP}^2$ at $n(4\leq n\leq
8)$ points. At this time $\langle K_{\widetilde{S}}\rangle^{\perp}$
is a root system of $E_n$-type. Thus $S$ is a blow-up of
$\mathbb{CP}^2$ at $n+1(4\leq n\leq 8)$ points.\par
   Secondly, suppose $C^2=0$. Then by Riemann-Roch theorem, the linear system $|C|$ defines a
   ruling over $\mathbb{P}^1$ with fiber $C$. Contract all $(-1)$ curves in fiber, we
   obtain a relatively minimal model (not unique), which is $\mathbb{P}^1\times
   \mathbb{P}^1$ or the Hirzebruch surface $\mathbb{F}_1$. So, $S$
   is a blow-up of $\mathbb{P}^1\times \mathbb{P}^1$ or
   $\mathbb{F}_1$ at $n$ points. And the lattice $\langle K_S,C\rangle^{\perp}$
   must be of $D_n$-type by Lemma~\ref{root-D}.\par
   Thirdly, suppose $C^2=1$. Then blow up one point $p_0\in C$, we obtain $\widetilde{S}$ which is a
   ruling over $\mathbb{P}^1$ with fiber $\widetilde{C}=C-E$ and section
   $E$ where $E$ is the exceptional curve associated to this
   blow-up. Contracting all $(-1)$ curves in fiber which do not
   intersect with $E$, we will obtain $\mathbb{F}_1$. Thus
   $\widetilde{S}$ is a blow-up of $\mathbb{F}_1$ at $n$ points. And we
   have
   $\langle K_{S},C\rangle^{\perp}\cong
   \langle K_{\widetilde{S}},\widetilde{C},E\rangle^{\perp}$. Therefore the
   lattice is a root lattice of $A_n$-type by Lemma~\ref{root-A}.\par
    Finally, suppose $C^2\geq 2$. Note that since we assume $C^2\neq 4$, the situation of Lemma~\ref{root-A} (iv) (a)
    can not happen. So we only need to discuss the case where $C^2=
   2$, because the discussion on general cases is similar.
   Blowing up $S$ at two points $p,q\in C,p\neq q$, we obtain
   $\widetilde{S}$ with exceptional curves $E_p,E_q$. Let
   $\widetilde{C}=C-E_p-E_q$ be the strict transform of $C$, then
   $|\widetilde{C}|$ defines a ruling with fiber $\widetilde{C}$ and section $s=E_p$ (fixed).
   Similarly as before, contracting all $(-1)$ curves $E$ in fiber which
   satisfy $E\cdot \widetilde{C}=0=E\cdot s$, we will obtain
   $\mathbb{F}_1$. Then $\widetilde{S}$ can be considered as a
   blow-up of $\mathbb{F}_1$ at $n$ points.  Note that $\langle K_{S},C\rangle^{\perp}\cong
   \langle K_{\widetilde{S}},\widetilde{C},s,E_q\rangle^{\perp}$. We know that
   $\langle K_{\widetilde{S}},\widetilde{C},s\rangle^{\perp}$ is a root lattice
   of $A_n$-type from Lemma~\ref{root-A}. Then the result follows also from Lemma~\ref{root-A}. $\hfill\Box$

\begin{remark} {\rm We extend the definition of $E_n$ surfaces to all $n$ with
$0\leq n\leq 8$, by defining $E_n(n\leq 3)$ surfaces to be del Pezzo
surfaces of degree $9-n$.}\end{remark}

 \begin{corollary} \label{excep-sys-effect} On an $ADE$ surface, any
exceptional divisor perpendicular to $C$ is represented by an
irreducible curve. Therefore, any exceptional system consists of
exceptional curves.\end{corollary}

\noindent{\bf Proof}. In $E_n$ case, the result follows from
Proposition~\ref{ADE-class} and Lemma~\ref{equiv-del-pezzo}. In $D_n$ and $A_n$ cases, according to
Proposition~\ref{ADE-class}, the result is obvious. $\hfill\Box$\\

    In the following we generalize the definition for $n\leq 8$ points being {\it in general position} to any $n\geq 0$.
Denote $S=\mathbb{P}^2$ (or $\mathbb{P}^1\times \mathbb{P}^1$ or
$\mathbb{F}_1$). Denote $S_n(x_1,\cdots,x_n)$ (or $S_n$ for brevity)
the blow-up of $S$ at $n$ points $x_1,\cdots,x_n$. We say that
$x_1,\cdots,x_n$ are {\it in general position} if any smooth
rational curve on $S_n$ has a self-intersection number at least
$-1$. And we say that $x_1,\cdots,x_n$ are {\it in almost general
position} if any smooth rational curve on $S_n$ has a
self-intersection number at least $-2$.

\begin{corollary} Let $(S,C)$ be an $ADE$ surface.\par
   (i) In $E_n$ case, blowing down the $(-1)$ curve $C$, we obtain a
del Pezzo surface of degree $9-n$. \par (ii) In $D_n$ case, $S$ is
just a blow-up of  $\bbP^1\times \bbP^1$ or $\bbF_1$ at $n$ points
in general position with $C$ as the natural ruling. \par (iii) In
$A_n$ case, let $\widetilde{S}$ be the blow-up of $S$ at a point on
$C$, with the exceptional curve $E$, then $\widetilde{{S}}$ is a
blow-up of $\bbF_1$ at $n+1$ points, and the strict transform
$\widetilde{C}$ of $C$ defines a ruling with $E$ as the section of
$\bbF_1$. $\hfill\Box$\end{corollary}

\section{Lie algebra bundles over rational surfaces of $ADE$-type and their representation bundles}

 When $G$ is of $ADE$-type, to each $ADE$ surface
 $S$, we can construct a natural
$\mathcal{G}=Lie(G)$ bundle and natural fundamental representation
bundles over $S$, which are determined by the lines (or exceptional
divisors in general) and rulings on $S$.

\begin{definition} \label{bundle} By a $Lie$ $algebra$ $\mathcal{G}=Lie(G)$
$bundle$, we mean a vector bundle which fiberwise carries a Lie
algebra structure of $\mathcal{G}$-type, and this Lie algebra
structure is compatible with trivialization of
 this bundle. By a $representation\  bundle$ of a $\mathcal{G}$ bundle,
 we mean a vector bundle $\mathcal{V}$ which fiberwise is a representation of
 $\mathcal{G}$, and the action of $\mathcal{G}$ on $\mathcal{V}$ is
 compatible with their trivialization.\end{definition}

   We describe these bundles in the following, and give the detailed arguments just in $E_n$ case,
since other cases are similar.

\subsection{$E_n$ bundles over $E_n$ surfaces}

    Let $(S,C)$ be an $E_n$ surface. Recall that $S=X_{n+1}(x_1,\cdots,x_{n+1})$
   where  $C$ be the exceptional divisor associated to the
   blow-up at $x_{n+1}$. Denote
   $\widetilde{S}=X_{n}(x_1,\cdots,x_n)$. Since $\langle K_{S},C\rangle^{\perp}\cong
   K_{\widetilde{S}}^{\perp}$, we can just consider the surface
   $\widetilde{S}=X_{n}(x_1,\cdots,x_n)$.

   Since we have a root system of $E_n$-type attached to $X_n$,
   inspired by the Cartan decomposition of a complex simple Lie algebra,
   we can construct a Lie algebra bundle over $X_n$ as follows:
       $$\mathscr{E}_n=\mathcal{O}^{\oplus n}\bigoplus_{D\in R_n}\mathcal{O}(D).$$

   The fiberwise Lie algebra structure of $\mathscr{E}_n$ is defined as the following.
   Fix the system of simple roots of $R_n$ as
   $$\Delta(E_n)=\{\alpha_1=l_1-l_2,\alpha_2=l_2-l_3,\alpha_3=h-l_1-l_2-l_3,\cdots,\alpha_n=l_{n-1}-l_n \},$$
   and take a trivialization  of $\mathscr{E}_n$. Then over a trivializing open subset $U$,
   $\mathscr{E}_n|_U\cong U\times(\mathbb{C}^{\oplus n}\bigoplus_{\alpha\in
   R_n}
   \mathbb{C}_{\alpha})$. Take a Chevalley basis $\{x_{\alpha}^U,\alpha\in R_n;h_i,1\leq i\leq n\}$
   for $\mathscr{E}_n|_U$  and define the Lie algebra structure by the  following
   four relations, namely, Serre's relations on Chevalley basis (see \cite{Humph1},
   p147):\\

 (a) $[h_ih_j]=0,1\leq i, j\leq n;$\par
 (b) $[h_ix_{\alpha}^U]=\langle \alpha,\alpha_i\rangle x_{\alpha}^U,1\leq i\leq n,\alpha\in
     R_n;$\par
 (c) $[x_{\alpha}^Ux_{-\alpha}^U]=h_{\alpha}$  is a $\mathbb{Z}$-linear combination of
     $h_1,\cdots,h_n;$\par
 (d) If $\alpha,\beta$ are independent roots, and
     $\beta-r\alpha,\cdots,\beta+q\alpha$ are the
     $\alpha$-string through $\beta$,
     then $[x_{\alpha}^Ux_{\beta}^U]=0$ if $q=0$, while $[x_{\alpha}^Ux_{\beta}^U]=\pm(r+1)x_{\alpha+\beta}^U$ if
     $\alpha+\beta\in R_n.$\\

   Note that $h_i,1\leq i\leq n$ are independent of any trivialization,
so the relation $(a)$ is always invariant under different
trivializations. If $\mathscr{E}_n|_V\cong
V\times(\mathbb{C}^{\oplus n}\bigoplus_{\alpha\in R_n})$ is another
trivialization, and $f_{\alpha}^{UV}$ is the transition function for
the line bundle $\mathcal{O}(\alpha)(\alpha\in R_n)$, that is,
$x_{\alpha}^U=f_{\alpha}^{UV}x_{\alpha}^V$, then the relation $(b)$
is
  $$[h_i(f_{\alpha}^{UV}x_{\alpha}^V)]=\langle \alpha,\alpha_i\rangle f_{\alpha}^{UV}x_{\alpha}^V,$$
   that is, $$[h_ix_{\alpha}^V]=\langle \alpha,\alpha_i\rangle x_{\alpha}^V.$$ So $(b)$
is also invariant. $(c)$ is also invariant since
$(f_{\alpha}^{UV})^{-1}$ is the transition function for
$\mathcal{O}(-\alpha)(\alpha\in R_n)$. Finally, $(d)$ is invariant
since $f_{\alpha}^{UV}f_{\beta}^{UV}$ is the transition function for
$\mathcal{O}(\alpha+\beta)(\alpha,\beta\in R_n)$.\par Therefore, the
Lie algebra structure is compatible with the trivialization. Hence
it is well-defined. In other words,  we can construct globally a Lie
algebra bundle over a
surface once we are given a root system consisting of divisors on this surface.\\

   The following relations  are intricate. One is the relation between $I_n$ (the set of all exceptional
   divisors) and the fundamental representation associated to the
   highest weight $\lambda_n$ which is dual to the simple root
   $\alpha_n$ (see Figure 1).  Another one is the relation between the set of
   rulings  and the fundamental representation associated to the
   highest weight $\lambda_1$ which is dual to the simple root
   $\alpha_1$ (Figure 1). We explain the relations in the
   following.\par
        Let  $\mathbb{L}_n$ be
   the fundamental representation with the highest weight $\lambda_n$.  Then we
   have:\\

   \begin{center}
   \begin{tabular}{|c|c|c|c|c|c|c|c|c|}
     \hline
     n& 1 & 2 & 3 & 4 & 5 & 6 & 7 & 8 \\
     \hline
     $dim\ \mathbb{L}_n$  & 1 & 3 & 6 & 10 & 16 & 27 & 56 & 248 \\
     \hline
     $| I_n|$  & 1 & 3 & 6 & 10 & 16 & 27 & 56 & 240 \\
     \hline
   \end{tabular}\\
   \end{center}

      \vspace{0.5cm}Denotes $Ru_n$ the set of all rulings on $X_n$. Let  $\mathbb{R}_n$ be
   the fundamental representation with the highest weight $ \lambda_1$.  Then we
   have:\\

   \begin{center}
   \begin{tabular}{|c|c|c|c|c|c|c|c|c|}
     \hline
     n& 1 & 2 & 3 & 4 & 5 & 6 & 7 & 8 \\
     \hline
     $dim\ \mathbb{R}_n$  & 1 & 2 & 3 & 5 & 10 & 27 & 133 & 3875 \\
     \hline
     $|Ru_n|$  & 1 & 2 & 3 & 5 & 10 & 27 & 126 & 2160 \\
     \hline
   \end{tabular}\\
   \end{center}

   \vspace{0.5cm} Inspired by these, we can construct a fundamental representation
   bundle $\mathscr{L}_n$ (respectively $\mathscr{R}_n$) using the exceptional divisors (respectively the rulings)
   on $X_n$ as follows.

   \begin{eqnarray} \mathscr{L}_n &=& \bigoplus_{l\in I_n}\mathcal {O}(l) \mbox{ when } n\leq 7,\nonumber
   \\ \mathscr{L}_8 &=& \bigoplus_{l\in I_8}\mathcal {O}(l)\oplus\mathcal {O}(-K)^{\oplus 8}.\nonumber
   \end{eqnarray}
   Respectively,
   \begin{eqnarray}\mathscr{R}_n &=& \bigoplus_{R\in Ru_n}\mathcal {O}(R) \mbox{ when } n\leq 6,\nonumber
   \\ \mathscr{R}_7 &=& \bigoplus_{R\in Ru_7}\mathcal {O}(R)\oplus\mathcal {O}(-K)^{\oplus 7}.\nonumber
   \end{eqnarray}

       The fiberwise action is defined naturally, which is in fact compatible with any
       trivialization.\par
      For example we consider the bundle $\mathscr{L}_n$ and suppose $n\leq 7$.
   Take $U,V$ as before, and suppose they also
   trivialize  $\mathscr{L}_n$, that is  $\mathscr{L}_n|_U \cong U\times(\bigoplus\limits_{l\in I_n}\mathbb{C}_l)$
   and $\mathscr{L}_n|_V \cong V\times(\bigoplus\limits_{l\in I_n}\mathbb{C}_l)$. Take $e_l^U$ (resp. $e_l^V=g^{VU}e_l^U$) to be
   the basis of $\mathbb{C}_l$ over $U$ (resp. $V$). Then define
   $x_{\alpha}^U.e_l^U$ to be equal to $e_{l'}^U$ if $l'=\alpha+l\in I_n$ and be equal to $0$ otherwise. And define
   $h_{\alpha}.e_l^U=(\alpha\cdot l)e_l^U$.\par
      Note that the situation here is slightly different from some standard usage,
   for example \cite{Bourbaki} \cite{Humph1}, since the self-intersection number
   of an element of $R_n$ or $I_n$ is negative. But this does not matter if we take the simple root system to be
   $\{-\alpha_1,\cdots,-\alpha_n\}$, and take the pairing to be $(x,y):=-(x\cdot y)$.
   Firstly since $ \lambda_n(-\alpha_i)=(-\alpha_i,l_n)=\alpha_i\cdot l_n=\delta_{in}$, we have $\lambda_n\cong(\cdot,l_n)$.
   Secondly the action is irreducible since the Weyl group acts on $I_n$ transitively.
   Lastly $e_{l_n}^U$ is the maximal vector of weight $\lambda_n$. Therefore this fiberwise
   action does define the highest weight module with the highest weight $\lambda_n$
   (see \cite{Humph1}). \par
     Obviously, this fiberwise  Lie algebra action is compatible with the trivialization. \\

      For $\mathscr{L}_8$, note that the bijection $I_8\rightarrow R_8$ given by
    $l\mapsto l+K$ induces an isomorphism $$\mathscr{E}_8\cong \mathscr{L}_8\otimes
    \mathcal{O}(K).$$ This implies $\mathscr{L}_8$ is just the adjoint
    representation bundle.\\

    Similarly, $\mathscr{R}_n$ is the fundamental representation
    bundle with the highest weight $\lambda_1\cong(\cdot,h-l_1)$ and the maximal vector  $e_{h-l_1}^U$,
    where the simple root system and the pairing are defined as
    above. We also have that $\mathscr{R}_7\otimes \calO(K)\cong \mathscr{E}_7$ is the adjoint representation bundle.

\begin{example} {\rm Let us look at the $sl(2)$ sub-bundle
$$\calO\oplus\calO(\alpha)\oplus(-\alpha),$$  where
$\alpha=l_1-l_2$. Then the bundle $\mathcal {O}(l_1)\oplus\mathcal
{O}(l_2)$ is the standard representation bundle. And the line bundle
$\mathcal {O}(h-l_1-l_2)$ is a trivial representation.}\end{example}

   In fact, the Lie algebra bundle $\mathscr{E}_n$ is uniquely
determined by its representation bundles $\mathscr{L}_n$ and
$\mathscr{R}_n$, according to \cite{Adams}. Concretely (see
\cite{Leung} for more details),
\par
  (i) $\mathscr{E}_4$ is the automorphism bundle of $\mathscr{R}_4$
  preserving $\wedge^5 \mathscr{R}_4\cong \mathcal{O}(-2K).$\par
  (ii) $\mathscr{E}_5$ is the automorphism bundle of $\mathscr{R}_5$
  preserving $q_5:\  \mathscr{R}_5\otimes\mathscr{R}_5\rightarrow \mathcal{O}(-K)$,
  where $q_5$ is defined by $\mathcal{O}(R')\otimes\mathcal{O}(R'')\rightarrow \mathcal{O}(-K)$
  if $R'+R''=-K$, and $0$ otherwise.\par
  (iii) $\mathscr{E}_6$ is the automorphism bundle of
  $\mathscr{R}_6$ and $\mathscr{L}_6$
  preserving

  $$\left \{ \begin{array}{l}c_6:\  \mathscr{L}_6\otimes\mathscr{L}_6\rightarrow
  \mathscr{R}_6,\ and \\c_6^{*}:\  \mathscr{R}_6\otimes\mathscr{R}_6\rightarrow
  \mathscr{L}_6\otimes \mathcal{O}(-K),
   \end{array}
\right. $$ where $c_6$ is defined by the map $(l_i,\ l_j)\mapsto
2h-\sum\limits_{k\neq i,j}l_k$ and $c_6^{*}$ is defined by the map
$(h-l_i,h-l_j)\mapsto h-l_i-l_j$.\par

   (iv)  $\mathscr{E}_7$ is the automorphism bundle of $\mathscr{L}_7$
  preserving $$f_7:\  \mathscr{L}_7\otimes\mathscr{L}_7\otimes\mathscr{L}_7\otimes\mathscr{L}_7\rightarrow \mathcal{O}(-2K),$$
  where $f_7$ is defined by the map $(C_1,C_2,C_3,C_4)\mapsto -2K$
  if $C_1+C_2+C_3+C_4=-2K$, $0$ otherwise.\par

   (v) $\mathscr{E}_8$ is the automorphism bundle of $\mathscr{L}_8$
  preserving $$\mathscr{L}_8\wedge\mathscr{L}_8\rightarrow \mathscr{L}_8\otimes\mathcal{O}(-K).$$

    For $X_6$, the bijection $Ru_6\rightarrow I_6$ defined by $R\mapsto
    -(R+K)$ induces an isomorphism $\mathscr{R}_6\cong \mathscr{L}_6^{*}\otimes
    \mathcal{O}(-K)$, which is consistent with the duality between
    $\mathbb{L}_6$ and $\mathbb{R}_6$ for the Lie group $E_6$.\par


\subsection{$D_n$ bundles over rational ruled surfaces}

      \ Let $(S,C)$ be a $D_n$ surface. By Proposition~\ref{ADE-class}, $S$
      dominates $\mathbb{F}_1$ or $\mathbb{F}_0(=\bbP^1\times\bbP^1)$ with ruling $C$. We
      can suppose that $S$ dominates $\mathbb{F}_1$ since for
      another  case the arguments is the same. Thus $S=Y_n(x_1,\cdots,x_n)$ is
      the blow-up of $\mathbb{F}_1$ at $n$ points
      $x_i,i=1,\cdots,n$, where for any $i$, $x_i$ does not lie on the section $s$.

      Since $R_n$ is a root system of type $D_n$, the Lie algebra bundle can be constructed as follows.
    $$\mathscr{D}_n=\mathcal{O}^{\oplus n}\bigoplus_{D\in R_n}\mathcal{O}(D).$$

    Recall that in $D_n$ case,
    \begin{eqnarray}I_n &=& \{C|\ C^2=C\cdot K=-1,C\cdot f=0\} \nonumber
   \\  &=& \{l_i,f-l_i\ |\ i=1,\cdots,n\}.\nonumber
   \end{eqnarray}
     The fundamental representation with the highest weight $\lambda_n$, where $\lambda_n$
    is the fundamental weight corresponding to $\alpha_n=l_{n-1}-l_n$, is
    $$\mathscr{W}_n=\bigoplus_{C\in I_n}\mathcal{O}(C).$$
       In fact, $\mathscr{W}_n$ is the standard representation bundle of
       $\mathscr{D}_n$. \par

     Note that there are $n$ singular fibers, and each singular fiber
     is of the form $l_i+l_i'$ where $l_i'=f-l_i,i=1,\cdots,n$. The relation
     $$\mathcal{O}(l_i)\otimes \mathcal{O}(l_i')=\mathcal{O}(f)$$
      implies we can define a non-degenerated fiberwise quadratic form
      $$q_n:\mathscr{W}_n\otimes \mathscr{W}_n \rightarrow \mathcal{O}(f).$$  The two spinor bundles are defined as
    $$\mathcal{S}_n^{+}=\bigoplus_{S^2=S\cdot K=-1, S\cdot f=1}\mathcal{O}(S)\mbox{ and }
    \mathcal{S}_n^{-}=\bigoplus_{T^2=-2,T\cdot K=0, T\cdot f=1}\mathcal{O}(T).$$
    Moreover, there are all kinds of structures on these
    representation bundles, for example, the Clifford
    multiplication:
    $$\mathcal{S}_n^+\otimes\mathscr{W}_{n}^*\rightarrow \mathcal{S}_n^-\mbox{ and }\mathcal{S}_n^-\otimes\mathscr{W}_{n}\rightarrow \mathcal{S}_n^+.$$
    When $n=2m-1$ is odd,  we have isomorphism $$(\mathcal{S}_n^+)^*\otimes\mathcal{O}_{Y_n}((m-4)f-K)\cong \mathcal{S}_n^-.$$
    When $n=2m$ is even,  we have isomorphisms
    \begin{eqnarray}(\mathcal{S}_n^+)^*\otimes\mathcal{O}_{Y_n}((m-3)f-K) & \cong & \mathcal{S}_n^+,\nonumber
   \\ (\mathcal{S}_n^-)^*\otimes\mathcal{O}_{Y_n}((m-4)f-K) & \cong & \mathcal{S}_n^-.\nonumber
   \end{eqnarray}
       For more details, see \cite{Leung}.

\subsection{$A_{n-1}$ bundles and their representation bundles }

   \ Let $S$ be an $A_{n-1}$ surface. By Proposition~\ref{ADE-class}, we can assume that
    $S=Z_n(x_1,\cdots,x_n)$ be  the blow-up of $\mathbb{F}_1$ at $n$ points
      $x_i,i=1,\cdots,n$, where for any $i$, $x_i$ does not lie on the section
      $s$. Recall that

      \begin{align}
            R_{n-1} &=\{l_i-l_j|\ i\neq j\} \mbox{ and } \nonumber\\
            I_{n-1} &=\{l_1,\cdots,l_n\}.\nonumber
      \end{align}
    Since $R_{n-1}$ is a root system of $A_{n-1}$-type, the Lie algebra bundle
    can be constructed as
$$\mathscr{A}_{n-1}=\mathcal{O}^{\oplus n-1}\bigoplus_{D\in
R_{n-1}}\mathcal{O}(D).$$

  And the standard representation bundle is
$$\mathcal{V}_{n-1}=\bigoplus_{C\in I_{n-1}}\mathcal{O}(C)=\bigoplus\limits_{i=1}^{n}\mathcal{O}(l_i).$$
  The $k^{th}$ fundamental representation bundle is just
  $$\wedge^k(\mathcal{V}_{n-1})\cong\bigoplus\limits_{i_1<\cdots< i_k}\calO(l_{i_1}+\cdots +l_{i_k}).$$ We also have
  $\mathscr{A}_{n-1}=\mathcal {E}nd_0(\mathcal{V}_{n-1})$.

    We summarize the content of this section as the following form.

\begin{conclusion} For every $ADE$ surface $S$, there is
     a  natural Lie algebra bundle of corresponding $ADE$-type over $S$.
     Furthermore, we can construct two natural fundamental representation
     bundles over $S$, using lines and rulings on $S$. Moreover, the Lie algebra bundle can be considered as
     the automorphism (Lie algebra) bundle of these fundamental representation
     bundles preserving natural structures. $\hfill\Box$\end{conclusion}

\section{Flat $G$ bundles over elliptic curves}

     In this section we review some well-known results about flat $G$ bundles over
elliptic curves.\par
     Let $\Sigma$ be an elliptic curve with identity element $0$.   The fundamental group
$\pi_1(\Sigma)=\mathbb{Z}\oplus\mathbb{Z}$.  Let $G$ be a compact,
simple and simply connected Lie group of rank $r$ with root system
$R$, coroot system $R_c$, Weyl group $W$, root lattice $\Lambda$,
coroot lattice $\Lambda_c$ and maximal torus $T$. The dual lattice
$\Lambda_c^{\vee}$ of $\Lambda_c$ is the weight lattice. We denote
the moduli space of flat $G$-bundles over $\Sigma$ by
$\mathcal{M}_{\Sigma}^{G}$. It is well-known that we have the
following  isomorphisms.
\begin{align}
 \mathcal{M}_{\Sigma}^{G}&\cong Hom(\pi_1(\Sigma),G)/ad(G) \nonumber\\
                         &\cong Hom(\pi_1(\Sigma),T)/W \nonumber\\
                         &\cong T\times T/W \nonumber\\
                         &\cong \Sigma\otimes_{\mathbb{Z}}\Lambda_c/W. \nonumber
\end{align}

The second isomorphism is because of Borel's theorem \cite{Borel}
which says that a commuting pair of elements in $G$ can be
diagonalized simultaneously. The last isomorphism comes from
$$Hom(\pi_1(\Sigma),T)=Hom(\pi_1(\Sigma),U(1)\otimes_{\mathbb{Z}}\Lambda_c)\cong
Hom(\pi_1(\Sigma),U(1))\otimes_{\mathbb{Z}}\Lambda_c$$
 and $$Hom(\pi_1(\Sigma),U(1))\cong Pic^0(\Sigma)\cong \Sigma.$$

 A theorem of Bernshtein-Shvartsman \cite{BS} and Looijenga \cite{Looijenga} says that
  $$\Sigma\otimes_{\mathbb{Z}}\Lambda_c/W\cong\mathbb{WP}^r_{s_0=1,s_1,\cdots,s_r},$$
  where the latter is the weighted projective space with weights
  $s_i$'s, and $s_1,\cdots, s_r$ are the coefficients of the highest
  coroot of $R_c$.\par
    One element of $Hom(\Lambda,\Sigma)/W$ can only determine a flat $ad(G)=G/C(G)$ bundle in
    general.
    For the adjoint group  $ad (G)$,   the moduli space of flat $ad (G)$ bundles $\mathcal{M}_{\Sigma}^{ad(G)}$
    contains
    $Hom(\Lambda,\Sigma)/W$ as a connected component (see \cite{FMW1}).
    On the other hand, we have the following short exact
sequences:
$$0\rightarrow \Lambda\rightarrow\Lambda_c^{\vee}\rightarrow \Gamma\rightarrow
0$$ and $$0\rightarrow Hom(\Gamma,\Sigma)\rightarrow
Hom(\Lambda_c^{\vee},\Sigma)\rightarrow Hom(\Lambda,\Sigma)
\rightarrow 0.$$ Here $\Gamma$ is a finite abelian group. The second
sequence is exact since $\Sigma$ is a divisible abelian group. It
follows that $Hom(\Lambda,\Sigma)$ and
$\Sigma\otimes_{\mathbb{Z}}\Lambda_c$ are isogenous as abelian
varieties. Let $d$ be the exponent of the finite group $\Gamma$. If
we fix a $d^{th}$ root of unity in $Jac(\Sigma)\cong \Sigma$ then we can
extend uniquely a homomorphism $f_0\in Hom(\Lambda,\Sigma)$ to a
homomorphism $f\in Hom(\Lambda_c^{\vee},\Sigma)\cong \Lambda_c\bigotimes\limits_{\bbZ}\Sigma$.
We have explained the following

\begin{lemma}\label{extend}
    When we fix a $d^{th}$ root of unity in $Jac(\Sigma)$, we have an isomorphism
    $$Hom(\Lambda_c^{\vee},\Sigma)/W\cong Hom(\Lambda,\Sigma)/W,$$ and therefore
$$\mathcal{M}_{\Sigma}^G\cong Hom(\Lambda,\Sigma)/W.$$
\end{lemma}

\begin{remark} {\rm We have constructed $ADE$ (Lie algebra) bundles over
$ADE$ rational surfaces.  We will see that the restriction of such a
Lie algebra bundle to the anti-canonical curve $\Sigma$ will
uniquely determine a flat $G$ bundle over $\Sigma$. To obtain a
simple Lie group $G=E_n$ (resp. $D_n$), we need to assume that
$4\leq n\leq 8$ (resp. $n\geq 3$).}\end{remark}

\section{Flat $G$ bundles over elliptic curves and rational surfaces: simply laced cases}

   From this section on, we fix our $ADE$ surface $S$ to be the rational surface $X_n(x_1,\cdots,x_{n})$,
   $Y_n(x_1,\cdots,x_n)$, or $Z_n(x_1,\cdots,x_n)$. For $X_n$, we assume  $n\leq 8$. \par

      Given any smooth elliptic curve $\Sigma$ with identity $0\in \Sigma$,
we assume that our surface $S$ contains $\Sigma$ as an
anti-canonical curve. For this aim, we first embed $\Sigma$ into
$\mathbb{P}^2$ as an anti-canonical curve, using the projective
embedding $\phi$ determined by the linear system $| 3(0)|$  where
$(0)$ is the divisor of the identity element of $\Sigma$, and assume
that all these blown up points $x_i\in\Sigma$ for $i=1,\cdots,n$,
and that $0,x_1,\cdots,x_n$ are in general position. Moreover, we
blow up $\mathbb{P}^2$ at $0$ to obtain the embedding of $\Sigma$
into $\mathbb{F}_1$ as an anti-canonical curve, and take the
exceptional curve $l_0$ as the section $s$ for the ruled surface
$\mathbb{F}_1$.

\begin{convention} In $Z_n$ case, it is well-known that in order
to obtain a flat $SU(n)$-bundle over $\Sigma$ we need one more
assumption:
$$\sum x_i=0\mbox{ in }\Sigma.$$\end{convention}

 We explain  how  the moduli space
$\mathcal{M}_{\Sigma}^{G}$ is related to the moduli space of
rational surfaces of the above types. Denote $\mathcal{S}(\Sigma,G)$
the moduli space of the pairs $(S,\Sigma)$, where $S$ is an $ADE$
rational surface of type being the same as that of $G$ and
$\Sigma\in|-K_S|$.

\begin{proposition} \label{moduli-map} There exists a well-defined map
$$\phi:\ \mathcal{S}(\Sigma,G)\rightarrow
Hom(\Lambda,\Sigma)/W,$$ where $\Lambda$ is the lattice $P_n$ or
$P_{n-1}$ defined in Section 1.\end{proposition}

\noindent{\bf Proof}. First we consider the case where $S=X_n$ is a
Del Pezzo surface, that is, all blown up points are in general
position. Suppose we are given the pair $(X_n,\Sigma\in  |-K_{X_n}|\
)$. For each element $y\in P_n$, $y$ stands for a holomorphic line
bundle over $S$. Restricting $y$ to $\Sigma$, we obtain a
holomorphic line bundle over $\Sigma$, denoted by $\mathcal {L}_y$.
The degree of $\mathcal {L}_y$ is
$$deg(\mathcal {L}_y)=y\cdot (-K)=0.$$ So $\mathcal {L}_y$ is an
element of the Jacobian of $\Sigma$, which is canonically isomorphic
to $\Sigma$ since the identity element of $\Sigma$ is given. Thus we
obtain a map from $P_n$ to $\Sigma: y\mapsto \mathcal {L}_y$, which
is obviously a homomorphism of abelian groups. But for one pair
$(X_n,\Sigma)$, we can have different choices of simple roots in
order to identify $P_n$ with the root lattice of $E_n$, and all
choices are only differed by the action of the Weyl group $W(E_n)$.
So finally we obtain a well-defined map from the moduli space
$\mathcal{S}(\Sigma,E_n)$ of such pairs $(X_n,\Sigma)$ to the
projective variety $Hom(P_n,\Sigma)/W(E_n)$.\par
 The other two cases are similar. Roughly speaking, given a pair $(Y_n,\Sigma)$  (resp. $(Z_n,\Sigma)$),
 we obtain an element in\\
  $$Hom(P_n,\Sigma)/W(D_n)\  (\mbox{ resp. }Hom(P_{n-1},\Sigma)/W(A_{n-1})).\eqno\Box $$\\

    In fact we can prove a theorem of Torelli type for the above correspondings. Roughly speaking,
 the moduli space of the pairs $(S,\Sigma)$ is isomorphic to
 $$Hom(\Lambda,\Sigma)/W,$$ where $\Lambda$ is our root
lattice.

\begin{definition} \label{config} Let $S=X_n$,  $Y_n$, or $Z_n$.  An exceptional
system $\zeta_n=(e_1,\cdots,e_n)\in C_n$ on $X_n$ (resp. $Y_n$,
$Z_n$) is called  a $G$-$configuration$ for $G=E_n$ (resp. $D_n$,
$A_{n-1}$) if $e_n$ is a $(-1)$ curve, and after blowing down $e_n$,
$e_{n-1}$ is a $(-1)$ curve. And this process can be proceeded
successively until after blowing down $e_1$, we obtain
$\mathbb{P}^2$ (resp. $\mathbb{F}_1$) for $G=E_n$ (resp. $D_n$ and
$A_{n-1}$). Denote $\zeta_{G}$ a $G$-$configuration$. When $S$ is
equipped with a $G$-configuration $\zeta_G$, and $S$ has $\Sigma$ as
an anti-canonical curve, we call $S$ a $rational$ $ surface$ $ with$
$G$-$configuration$ and denote it by a pair $(S,G)$.
\end{definition}

   Equivalently, a $G$-configuration $\zeta_{E_n}$ (resp. $\zeta_{D_n}$ or
   $\zeta_{A_{n-1}}$) on $S=X_n$ (resp. $Y_n$, $Z_n$), means that
  $S$ could be considered as the blow-up of $\mathbb{P}^2$
(resp. $\mathbb{F}_1$,  $\mathbb{F}_1$) at $n$ (maybe not distinct)
points $y_1,\cdots,y_n\in S$ successively, such that
$e_1,\cdots,e_n$ are the corresponding exceptional divisors.

\begin{lemma}\label{config-geq(-2)} Let $S$ be a surface with $G$-configuration. Then any smooth
rational curve on $S$ has a self-intersection number at least $-2$.
Furthermore, in $E_n$ case, all these $(-2)$ curves form chains of
$ADE$-type.\end{lemma}

 \Prf. Let $L$ be a smooth rational curve on
$S$. Then $L\cdot\Sigma\geq 0$. By adjoint formula, we have
$-2=L^2+L\cdot K_S$. Since $\Sigma$ is linearly equivalent to
$-K_S$, we have $L^2\geq -2$. For the last assertion, see
\cite{Demazure}. \hf \\

   On an $ADE$ surface, by Corollary~\ref{excep-sys-effect}, any exceptional system is an $ADE$-configuration.
  Thus, we can restate the result of Lemma~\ref{root-E} (ii), Lemma~\ref{root-D} (ii) and Lemma~\ref{root-A} (ii)
  as follows.

\begin{proposition} \label{W-action} For an $ADE$ surface, $W(G)$ acts on the set
of all $G$-configurations simply
 transitively.$\hfill\Box$\end{proposition}

   This proposition implies that a $G$-configuration determines exactly an isomorphism from $P_n$
(or $P_{n-1}$ for $A_{n-1}$) to the corresponding
root lattice $\Lambda(G)$.\\

   An $A_{n-1}$-configuration on $Z_n$ is illustrated in the following figure\\

  \begin{center}
  \includegraphics{ade.5}\\
  \end{center}

   A $D_n$-configuration on $Y_n$ is illustrated in the following figure\\

   \begin{center}
   \includegraphics{ade.6}\\
   \end{center}

  And an $E_n$-configuration on $X_n$ is illustrated in the following figure\\

   \begin{center}
   \includegraphics{ade.7}\\
   \end{center}

     Recall the definition for $\zeta_{D_n}$: $\zeta_{D_n}=(e_1,\cdots,e_n)$ where
     $e_i\cdot K_{Y_n}=-1,\ e_i\cdot f=0,\  e_i\cdot e_j=\delta_{ij}$ and $\sum e_i\cdot s\equiv 0 \ mod \ 2$. Next we explain
   geometrically why we need to assume that $\sum e_i\cdot s\equiv 0 \ mod \
   2$.

\begin{definition} Let $C\subset \mathbb{P}^2$ be a curve of
degree $d$. A point $P\in C$ is called a $ordinary$ $k$-$fold$
$point$  of $C$ if $P$ is a $k$-fold singular point and $C$ has $k$
distinct tangent directions at $P$. \end{definition}

\begin{lemma}\label{ordinary-sing-curve} Let $C$ be a plane curve of degree $d$ with an
ordinary $(d-1)$-fold point $P$. Then\par (i) $P$ is the only
singular point of $C$.\par (ii) The normalization of $C$ is a smooth
rational curve.\par (iii) Fix a point $P\in \mathbb{P}^2$. Then  the
variety of all plane curves of degree $d$ with $P$ as an ordinary
$(d-1)$-fold point is of dimension $2d$.\par (iv) Given $P$ and
other $2d$ generic points, there exists a unique curve $C\subset
\mathbb{P}^2$ of degree $d$, such that $C$ has $P$ as an ordinary
$(d-1)$-fold point and passes through these $2d$ generic points.
\end{lemma}

\noindent{\bf Proof}. (i) Apply Bezout's theorem. (ii) Apply the
genus formula. (iii) Let $[x,y,z]$ be the homogenous coordinates of
$\mathbb{P}^2$, and $P=[1,0,0]$. Then $C$ is defined by the
polynomial
$$f(x,y,z)=g(y,z)+\prod\limits_{i=1}^{d-1}(a_iy-b_iz)x,$$ where
$deg(g)=d$. Therefore, the dimension is $2d$. $\hfill\Box$

\begin{proposition}  Let $\Sigma$ be embedded into $\mathbb{F}_1$
(with section $s$) as a smooth anti-canonical curve and
$x_1,\cdots,x_n$ are distinct points of $\Sigma$. Blowing up
$\mathbb{F}_1$ at $x_i$'s we obtain $Y_{n}$ with corresponding
exceptional curves $l_i,i=1,\cdots,n$.
\par
  (i) When $n=2k$, if $x_1,\cdots,x_n$ are in general position, then
after  contracting $f-l_1,\cdots,f-l_n$, we still obtain the surface
$\mathbb{F}_1$. In other words, we obtain the same surface $Y_{2k}$
 by blowing up either $\{x_1,\cdots,x_n\}$, or $\{-x_1,\cdots,-x_n\}$. \par
  (ii) When $n=2k+1$, if $x_1,\cdots,x_n$ are in general position, then
after  contracting $f-l_1,\cdots,f-l_n$, the resulting surface is
$\mathbb{P}^1\times \mathbb{P}^1$, but not
 $\mathbb{F}_1$. \end{proposition}

\noindent{\bf Proof}. Let $C$ be a negative rational curve in
$Y_{n}$ which doesn't intersect $f-l_i,i=1,\cdots,n$. Then $C$
satisfies the following equations
$$\left \{ \begin{array}{l}C\cdot C = -m,m>0;
   \\ C\cdot K = m-2;
   \\ C\cdot (f-l_i) = 0,i=1,\cdots,n.
      \end{array}
\right.$$

 Since $C$ is a rational curve and $\Sigma\in|-K|$, $C\cdot
(-K)\geq 0$. So $m\leq 2$. Then $m=1\mbox{ or }2$. Considering
$\mathbb{F}_1$ as the blow-up of $\mathbb{P}^2$ at $0\in \Sigma$
with exceptional curve $s$, we can assume $C=a\cdot h-b\cdot s-\sum
c_i\cdot l_i,a\geq 0,b\geq 0,c_i\geq 0$. Solving the system of
equations, we obtain

$$\left \{ \begin{array}{l}m = 1\mbox{ or }2,
   \\b=a-1,
   \\ c_i = 1,i=1,\cdots,n,
   \\ a = (n-1+m)/2.
      \end{array}
\right.$$

    For $m=1$, $n=2a$ is even. The class
$$C=ah-(a-1)s-\sum\limits_{i=1}^{n=2a} l_i=af+s-\sum\limits_{i=1}^{2a} l_i.$$
    This means that all of the points $0,x_1,\cdots,x_n$ lie on
the curve $\pi(C)$, where $\pi: Y_n\rightarrow \mathbb{P}^2$ is the
blow-up of $\mathbb{P}^2$ successively at $0,x_1,\cdots,x_n$. There
exists exactly one such curve $C$ for generic $x_1,\cdots,x_n$, and
it is smooth, by Lemma~\ref{ordinary-sing-curve}. Hence, after
contracting $f-l_1,\cdots,f-l_{2a}$, we still obtain $\mathbb{F}_1$.

   For $m=2$, $n=2a+1$ is odd. The class
$$C=ah-(a-1)s-\sum\limits_{i=1}^{n=2a+1} l_i=af+s-\sum\limits_{i=1}^{2a+1} l_i.$$
    This means that all of the points $0,x_1,\cdots,x_n$ lie on
the curve $\pi(C)$, where $\pi: Y_n\rightarrow \mathbb{P}^2$ is the
blow-up of $\mathbb{P}^2$ successively at $0,x_1,\cdots,x_n$. There
exists no such curves for generic $x_1,\cdots,x_n$, by
Lemma~\ref{ordinary-sing-curve}.
  Hence,  after contracting $f-l_1,\cdots,f-l_{2a+1}$, no rational curves
  with negative self-intersection number can survive.
  Therefore the resulting surface  is $\mathbb{P}^1\times\mathbb{P}^1$, but not
 $\mathbb{F}_1$. $\hfill\Box$

\begin{example} {\rm Blowing up $\mathbb{F}_1$ at $2$ points $x_1,x_2$
we obtain $Y_2$. Contracting $f-l_1$ and $f-l_2$, or contracting
$l_1$ and $l_2$, we always obtain the surface $\mathbb{F}_1$. But
contracting $f-l_1$ and $l_2$, we just obtain the surface
$\mathbb{P}^1\times \mathbb{P}^1$ , but not $\mathbb{F}_1$!}
\end{example}

\begin{remark} {\rm (i) Lemma~\ref{ordinary-sing-curve} has a corresponding version for
$\mathbb{P}^1\times \mathbb{P}^1$.\par
    (ii)  A $G$-configuration $\zeta_G=(e_1,\cdots,e_n)$ for $S=X_n$ (resp.
$Y_n,Z_n$) just means that after blowing down
$e_n,e_{n-1},\cdots,e_1$ successively, we still obtain
$\mathbb{P}^2$ (resp. $\mathbb{F}_1$, $\mathbb{F}_1$).} \end{remark}

   Let $S$ be an $ADE$ surface equipped with a $G$-configuration $\zeta_G$.
we denote the moduli space  of the pairs $(S, \Sigma)$ by $\mathcal
{S}(\Sigma,G)$, where two pairs $(S,\Sigma)$ and
$(S^{'},\Sigma^{'})$ are equivalent if and only if there is an
isomorphism $\pi$ from $S$ to $S'$ such that $\pi|_{\Sigma}$  is
also an isomorphism from $\Sigma$ to $\Sigma'$.\par

  We show that $\mathcal {S}(\Sigma,G)$ is isomorphic to an
  open dense subset $U$ of the variety $Hom(\Lambda, \Sigma)/W$.
  In fact, for any element $\theta\in (Hom(\Lambda,\Sigma)/W)\backslash U$, the boundary
  component,  we can find possibly
non-equivalent pairs $(S,\Sigma)$ such that $\theta$ comes from the
restriction. Thus, we can complete $\mathcal {S}(\Sigma,G)$ by
adding these pairs and identifying them as one point. Denote the
completion by $\overline{\mathcal {S}(\Sigma,G)}$. Then we can
identify $\overline{\mathcal {S}(\Sigma,G)}$ with the projective
variety $Hom(\Lambda, \Sigma)/W$. This provides a natural
compactification for the moduli space $\mathcal {S}(\Sigma,G)$.

   More precisely, let $S=X_n$ (respectively, $Y_n$, $Z_n$) be an $ADE$ surface and $\Lambda$ be
the root lattice of $E_n$ (respectively, $D_n$, $A_{n-1}$) with
corresponding Weyl group $W$. And we fix a $3^{rd}$ (respectively,
$2^{nd}$, $n^{th}$) root of unity in $Jac(\Sigma)\cong \Sigma$
 in $E_n$ (respectively, $D_n$, $A_{n-1}$) case.  Then we have

 \begin{theorem} \label{main} (i) There is
an injective map $\phi$ from the moduli space $\calS(\Sigma,G)$ onto
an open dense subset of $Hom(\Lambda, \Sigma)/W$.\par
    (ii) $\phi$ can be extended to a bijective map
from the completion $\overline{\mathcal {S}(\Sigma,G)}$ onto
$Hom(\Lambda, \Sigma)/W$. \par
   (iii) Moreover, the completion is obtained by including  all rational surfaces with
$G$-configurations to $\mathcal {S}(\Sigma,G)$. Any smooth rational
curve on a surface corresponding to a boundary point has a
self-intersection number at least $-2$, and in $E_n$ case these
$(-2)$ curves form chains of $ADE$-type. \end{theorem}

\noindent{\bf Proof}. First we suppose $S=X_n$. We have constructed the map $\phi$
in Proposition~\ref{moduli-map}. We prove the injectivity. Fix a
$G$-configuration $\zeta_G=(l_1,\cdots,l_n)$ on $X_n$, and a simple
root system
$\alpha_1=l_1-l_2,\alpha_2=l_2-l_3,\alpha_3=h-l_1-l_2-l_3,\alpha_4=l_3-l_4,\cdots,\alpha_n=l_{n-1}-l_n$.
Blowing down $l_n,l_{n-1},\cdots,l_1$ successively, we obtain
$\mathbb{P}^2$ with $\Sigma$ as an anti-canonical curve. For all
$i=1,\cdots,n$, let $x_i\in X_n$ be the unique intersection points
of $l_i$ and $\Sigma$. Then $X_n$ can be considered as a blow-up of
$\mathbb{P}^2$  at these $n$ points $x_i\in \Sigma,i=1,\cdots,n$
with exceptional curves $l_i,i=1,\cdots,n$.

  According to previous arguments, we have a homomorphism $g\in
Hom(\Lambda, \Sigma)$. Let $g(\alpha_i)=p_i\in \Sigma$, then we have
the following equations by the group law of $\Sigma$ as an abelian
group
$$\left\{
\begin{array}{l}  x_1-x_2=p_1, \\x_2-x_3=p_2,
   \\-x_1-x_2-x_3=p_3,
   \\x_{k-1}-x_k=p_k,k=4,\cdots,n.
\end{array}\right.
$$
The determinant of the coefficient matrix of this system of linear
equations is $\pm 3$. So it has unique solution (if we fix a
$3^{rd}$ root of unity in $Jac(\Sigma$)). That is, $x_i$'s are uniquely determined
by $g$ up to Weyl group actions. The Weyl group actions just lead to
choices of other $G$-configurations. By Proposition~\ref{W-action},
this doesn't change the pair $(X_n,\Sigma)$. Hence, $\phi$ is
injective. These points $x_i$'s are not $''$in general position$''$
if and only if $p_i$'s will satisfy some (finitely many) equations.
That means the image of $\phi$ must be open dense in $Hom(\Lambda,
\Sigma)/W$. The extendability of $\phi$ is also because of the
existence and uniqueness of the solution of the above equations.\par

     For the cases of $Y_n$ and $Z_n$,
the arguments are similar. It is easy to see that the map $\phi$ is
well defined in both cases. For $Y_n$, the system of linear
equations is
$$\left \{ \begin{array}{l}-x_1-x_2=p_1,\\
      x_{k-1}-x_k=p_k,k=2,\cdots,n.
      \end{array}
\right.$$ The determinant is $\pm 2$. So the solution is uniquely
determined (if we fix a $2^{nd}$ root of unity in $Jac(\Sigma$)). The remained
arguments is just like the first case. At last, for the case of
$Z_n$, the system of equations is

$$\left \{ \begin{array}{l}\sum x_i=0,\\
      x_{k-1}-x_k=p_{k-1},k=2,\cdots,n.
      \end{array}
\right.$$

The determinant is $\pm n$. Then the solution is uniquely determined
(if we fix an $n^{th}$ root of unity in $Jac(\Sigma$)). The remaining arguments
are just the same as that in the $E_n$ case. These prove (i) and
(ii).
\par As for (iii), the
result follows from Lemma~\ref{config-geq(-2)}. \hfill$\Box$\\

\begin{remark}\label{deformation} {\rm The referee remarked that the set
$\phi(\mathcal{S}(\Sigma,G))$ in Theorem~\ref{main} was exactly the
complement of the discriminant in $Hom(\Lambda,\Sigma)/W$. This is
the case for $E_n$ type. As the referee indicated to us, this
follows from the description by Looijenga
\cite{Looijenga1}\cite{Looijenga2} and Pinkham \cite{Pinkham} of
$Hom(\Lambda,\Sigma)/W$ as the semi-universal deformation space of a
simple-elliptic singularity. The deformation space is realized as a
family of affine surfaces, and the fiberwise  compactification is a
Del Pezzo surface with an anticanonical elliptic curve as the
complement divisor. And the $-2$ curves on fibers produce the
vanishing cycles which determine the discriminant locus in
$Hom(\Lambda,\Sigma)/W$. For other cases, it is hoped to be true.
However, we can not give a proof at present. When the anticanonical
curve $C\in |-K_S|$ is a nodal rational curve, the moduli space of
pairs $(S,C)$ is considered by Looijenga in \cite{Looijenga98}. This
is in fact a degeneration of the situation above, where the elliptic
curve degenerates into a nodal curve. It is also interesting to
study the configurations  on such surfaces which are related to some
fundamental representations.}
\end{remark}

   As a conclusion of Lemma~\ref{extend} and Theorem~\ref{main}, we have

\begin{theorem}
  When we fix a $d^{th}$ root of unity in $Jac(\Sigma)$, we have a bijection
$$\overline{\mathcal
{S}(\Sigma,G)}\stackrel{\thicksim}{\longrightarrow}
\mathcal{M}_{\Sigma}^{G},$$ where $d$ is the exponent of the finite group $\Lambda_c/\Lambda$. \hfill$\Box$
\end{theorem}

\begin{remark} {\rm \cite{Tu}\cite{FMW1}\cite{FMW2}. The moduli space of flat
$A_n$ bundles over $\Sigma$ is exactly the ordinary projective space
$\mathbb{CP}^n$. This can be described as follows: a flat $SU(n+1)$
bundle is determined uniquely by $n+1$ points on $\Sigma$ with sum
equal to $0$, up to isomorphism. And $n+1$ points on $\Sigma$ with
sum equal to $0$ are determined uniquely by a global section
$H^0(\Sigma,\mathcal{O}_{\Sigma}(n(0)))$ up to scalar, where $(0)$ is the divisor of
the identity element $0$. So the moduli space of flat $SU(n+1)$
bundles is isomorphic to
$\mathbb{P}(H^0(\Sigma,\mathcal{O}_{\Sigma}((n+1)P)))=\mathbb{P}^n$.
From this we see that the moduli space of pairs $(S,\ \Sigma)$ is
just the ordinary complex projective space
$\mathbb{C}\mathbb{P}^n$.}\end{remark}

\begin{example} {\rm Let us look at what the pre-image of a trivial $G$-bundle is.
For example, in $E_8$ case, the trivial bundle means the element $0\in
Hom (\Lambda(E_8),\Sigma)/W(G)$. By the above correspondence,
 all $x_i=0$ in $\Sigma$. This means that we can blow up $\mathbb{P}^2$ at
 the identity element
 $0$ (an inflection point) eight times to obtain the surface represented by this
 pre-image, which is a boundary point in the moduli space $\overline{\mathcal {S}(\Sigma,G)}$.
 Blowing up once more, we obtain an elliptic fibration with
 a  singular fiber of $\widetilde{E_8}$-type \cite{BHPV}.}\end{example}

\end{document}